\documentclass[aap,MSNbibl,dvips]{arximspdf}
\usepackage{mathbh}
\usepackage{breakurl}

\doi{10.1214/10-AAP678}
\volume{20}
\issue{6}
\pubyear{2010}
\firstpage{1967}
\lastpage{1988}

\makeatletter
\def\R{\mathbb{R}}
\def\N{\mathbb{N}}
\def\P{\mathbb{P}}
\def\Z{\mathbb{Z}}

\def\e{\mathrm{e}}

\newcommand{\eqref}[1]{(\ref{#1})}
\newtheorem{theorem}{Theorem}

\newtheorem{proposition}{Proposition}
\newtheorem{lemma}{Lemma}
\newtheorem{corollary}{Corollary}
\newproclaim{remark}{Remark}
\makeatother

\begin{document}
\begin{frontmatter}

\title{Asymptotic regimes for the partition into colonies of a
branching process with emigration}
\runtitle{Partition into colonies}

\begin{aug}
\author{\fnms{Jean} \snm{Bertoin}\ead[label=e1]{jean.bertoin@upmc.fr}\thanksref{t1}}
\thankstext{t1}{Supported by ANR-08-BLAN-0220-01.}
\runauthor{J. Bertoin}
\affiliation{Universit\'e Pierre et Marie Curie}
\address{Laboratoire de Probabilit\'{e}s et\\
\quad Mod\`{e}les Al\'{e}atoires\\
Universit\'{e} Pierre et Marie Curie\\
4, Place Jussieu\\
F-75252 Paris Cedex 05\\
France\\
\printead{e1}} 
\end{aug}

\received{\smonth{8} \syear{2009}}
\revised{\smonth{1} \syear{2010}}

%
\begin{abstract}
We consider a spatial branching process with emigration in which
children either remain at the same site as their parents or migrate to
new locations
and then found their own colonies. We are interested in asymptotics of
the partition of the total population into colonies for large
populations with rare migrations.
Under appropriate regimes, we establish weak convergence of the
rescaled partition to some random measure that is constructed from the
restriction of a Poisson point measure to a certain random region, and
whose cumulant solves a simple integral equation.
\end{abstract}

%

\end{frontmatter}

%
\section{Introduction}\label{sec1}

Imagine a spatial branching process in which the child of an individual
either is a homebody, that is, remains at the same site as its parent, or
migrates to a new location which has never been occupied before and
then founds its own colony. We assume that the reproduction law is the
same for homebodies and migrants and do not depend on the spatial
location either, so this is essentially a discrete version of the
Virgin Island Model of Hutzenthaler \cite{H} when local competition
between individuals is discarded; see also \cite{HW} and references therein.

The dynamics of the process are entirely determined by
the pair $(\xi^{\mathrm{h}},\xi^{\mathrm{m}})$ of integer valued
random variables
giving the number of homebody children and the number of migrant
children of a typical individual.
The special case where each child chooses to emigrate with a fixed
probability $p\in(0,1)$ and independently of the other children can be
interpreted in the framework of the infinite-sites model in population genetics
by identifying a spatial location with a locus on a chromosome and $p$
with the rate of neutral mutations. This setting has motivated a number
of works in the literature; see in
particular \cite{GP} and \cite{Taib}.
In a quite different direction, we may also consider for instance the
cut-off situation where there is a threshold~$k$ such that the first
$k$ children of an individual are always homebodies while the next
children (if any) are forced to migrate. We may think of many other
simple rules as there are no assumptions on the correlation between
$\xi
^{\mathrm{h}}$ and $\xi^{\mathrm{m}}$.

We are interested in statistics of the decomposition of the entire
population according to the locations of individuals, which we call the
partition into colonies. This partition is naturally endowed with a
genealogical tree structure which has been described in \cite{Be1}.
Recent work \cite{Be2} focused on the special case where neutral
mutations occur in a Galton--Watson process with a fixed reproduction
law which is critical and has a finite variance. In asymptotic regimes
where the population is large and the mutation rate small, we
established a weak limit theorem for the tree of alleles, that is, in
the present setting, the partition into colonies equipped with its
genealogical structure.
The limit was described in terms of the genealogical tree of a
continuous state branching process in discrete time (cf. \cite{Jirina})
with an inverse Gaussian reproduction measure. In a related direction,
we also point at recent work by Abraham and Delmas \cite{ADV} on
pruning L\'evy continuum random trees.

In the present paper, we shall investigate more generally asymptotics
of the partition into colonies (ignoring its genealogical structure)
for branching processes with emigration when populations are large and
migrations rare. The regimes of interest are related to the well-known
limit theorems for rescaled Galton--Watson processes toward continuous
state branching processes in continuous time. Our main result (Theorem
\ref{T2}) states that after an appropriate rescaling, the partition
into colonies converges weakly to some random point measure. The latter
is constructed from the restriction of a Poisson point measure to a
certain random region. An important step in our analysis is that,
although in general the cumulant of this limiting random measure is not
explicitly known, it can be characterized as the unique solution to a
rather simple integral equation.

Let us briefly present the plan of this work by explaining our
approach. In Section \ref{sec2}, we point at the fact that the cumulant of the
partition into colonies solves a certain integral equation. This
equation stems from the extended branching property that is fulfilled
by the partition, and
is given in terms of the distribution of a pair of random variables
which arise naturally in this setting. We also recall a useful identity
in law which relates the preceding variables to that of passage times
in certain random walks.

In Section \ref{sec3}, we consider a sequence of branching processes with
emigration and introduce the basic assumptions. These are closely
related to the classical limit theorems for rescaled Galton--Watson
processes and involve L\'evy processes with no negative jumps.
Motivated by Section \ref{sec2}, we investigate limits in distribution for
passage times of random walks, and point at the role of the L\'evy
measure of a bivariate subordinator which arises in this setting.

In Section \ref{sec4}, we introduce a family of random point measures which are
constructed from Poisson point measures on a product space by
restriction to certain random domains. The main feature is that the
cumulant of such a random measure can be characterized as the unique
solution to another integral equation involving the intensity measure
of the underlying Poisson measure.

Our main result for limits in law of partitions into colonies is
presented and proved in Section \ref{sec5}. Roughly, we show that the cumulants
of the partitions into colonies of a sequence of branching processes
with emigration\vadjust{\goodbreak} converge after an appropriate rescaling to the unique
solution of an equation of the type which appeared in Section \ref{sec4}. More
precisely, it corresponds to the case where the intensity of the
driving Poisson measure is given by the L\'evy measure that has arisen
in Section \ref{sec3}.

Finally, Section \ref{sec6} is devoted to a few (hopefully) interesting
examples, partly to demonstrate the variety of possible asymptotic
behaviors. Roughly speaking, the common feature in these examples is
that the Galton--Watson process for which spatial locations of
individuals are ignored has a fixed distribution. We shall consider
different natural possibilities for selecting migrants children amongst
the progeny of an individual, which will yield different limiting
partitions into colonies. In the case corresponding to rare neutral
mutations in the infinite alleles branching process, the limiting
random partition
can be described in terms of certain Poisson--Kingman partitions which
have been considered by Pitman \cite{PiPK}.

\section{Preliminaries on partitions into colonies}\label{sec2}
In this section, we briefly introduce notation and present some basic
properties for Galton--Watson processes with emigration and the induced
partitions into colonies.
The material is essentially adapted from \cite{Be1} and \cite{Be2} to
which we refer for details, with the exception of Lemma \ref{L1} which
is new.

Roughly speaking, we consider a spatial haploid population model with discrete
nonoverlapping generations where each individual begets independently
of the others, according to a fixed reproduction law which is
independent of the location of that individual. We do not specify
geometrical details of the
space where individuals live as this would be irrelevant for the study;
the only implicit assumption is that this space is infinite. A child
can either stay at the same site as its parent or migrate to a new site
which has never been occupied before and then found its own colony.
This child is called a \textit{homebody} in the first case, and a {\it
migrant} in the second. For the sake of simplicity, we shall assume in
this work that at the initial time each ancestor lives in a different
location, although arbitrary initial conditions could be dealt with
more generally. The law of this model is thus entirely determined by
the number of ancestors and a pair of integer-valued random variables
$(\xi^{\mathrm{h}},\xi^{\mathrm{m}})$ which should be thought of as the
number of homebody children and the number of migrant children of a
typical individual.
For every $a\in\N$, we use the notation $\P_a$ for the probability
measure under which this model starts from $a$ ancestors.

If spatial locations are discarded, then the total number of
individuals per generation clearly forms a standard Galton--Watson
process with reproduction law given by the distribution of
$\xi=\xi^{\mathrm{h}}+\xi^{\mathrm{m}}$. We always assume that this
Galton--Watson process is critical or sub-critical, namely,
$\mathbb{E}(\xi)\leq1$, and implicitly exclude the degenerate case
where $\xi
\equiv1$, so the population becomes eventually extinct a.s.
The main object of interest in this work is the \textit{partition into
colonies}, which we represent as a random discrete measure
\[
{\mathcal P} = \sum_{j=1}^{\gamma} \delta_{C_j}.
\]
Here $\gamma$ is the total number of colonies (that is occupied sites)
and $C_j$ denotes the total number of individuals that lived at the
$j$th colony.
Observe that the first moment of ${\mathcal P}$ coincides with the
total population of the Galton--Watson process, namely,
\[
\zeta:=\sum_{j=1}^{\gamma} C_j = a + \sum_{k=1}^{\zeta} \xi_k,
\]
where $\xi_k=\xi^\mathrm{h}_k+\xi^\mathrm{m}_k$ stands for the
number of
children of the $k$th individual for some enumeration procedure,
and that the mass of ${\mathcal P}$ is just the number of colonies
\[
\gamma= a + \sum_{k=1}^{\zeta} \xi^{\mathrm{m}}_k.
\]

We denote the cumulant of partition into colonies when there is a
single ancestor
by
\[
K(f)=-\ln\mathbb{E}_1(\exp(-\langle{\mathcal P},f\rangle)),
\]
where $f\dvtx \N\to\R_+$ stands for a generic function and
\[
\langle{\mathcal P},f\rangle=\sum_{j=1}^{\gamma}f(C_j).
\]
We also point out from the branching property that for an arbitrary
number of ancestors $a\in\N$
we have
\[
\mathbb{E}_a(\exp(-\langle{\mathcal P},f\rangle))= \exp(-aK(f)),
\]
hence the cumulant $K$ characterizes the law of ${\mathcal P}$ under $\P
_a$ for any $a\geq1$.

The starting point of our analysis relies on the fact that this
cumulant is determined
in terms of the distribution of a pair of random variables
which appear naturally in the branching process with emigration.
Specifically, imagine for a while a variation of the model starting
from a single ancestor in which migrants are sterilized (i.e., they have
no offspring). We denote by $C$ the total number of individuals that
lived at the same site as the ancestor and by $M$ the number of
sterilized migrant children. In other words, $C$ is the size of the
colony generated by the ancestor and $M$ the number of colonies which
have been founded by migrant children of the ancestral colony.

\begin{lemma}\label{L1} For every function $f\dvtx \N\to\R_+$, the cumulant
$K(f)$ of ${\mathcal P}$ is the unique solution $\lambda\geq0$ to the equation
\[
\e^{-\lambda}=\mathbb{E}_1\bigl(\exp\bigl(-f(C)-\lambda M\bigr)\bigr).
\]
\end{lemma}

\begin{pf}
This stems from
the branching property which is inherited by the partition into
colonies. More precisely, we work under $\P_1$ and decompose the total
population
into the ancestral colony
and families generated by the migrant children of that colony.
Because the descent of each migrant child has the same distribution as
the initial spatial
Galton--Watson process, independently of the other migrant children
and of the homebody offspring of the ancestor, this yields
\begin{eqnarray*}
\exp(-K(f))&=&\mathbb{E}_1(\exp(-\langle{\mathcal P},f\rangle))\\
&=& \mathbb{E}_1\bigl(\exp(-f(C))\mathbb{E}_M(\exp(-\langle{\mathcal
P},f\rangle
))\bigr)\\
&=& \mathbb{E}_1\bigl(\exp\bigl(-f(C)-K(f)M\bigr)\bigr).
\end{eqnarray*}
We refer to Chauvin \cite{Chauvin} for a rigorous formulation of the
extended branching property of Galton--Watson processes at stopping
lines that we have used above, and also to \cite{Be1} for an
alternative argument based on the strong Markov property of random walks.

Uniqueness of the solution follows from the following observation.
Suppose first that
$f(C)\not\equiv0$. By H\" older's inequality, the map
\[
\lambda\to\lambda+ \ln\mathbb{E}_1\bigl(\exp\bigl(-f(C)-\lambda M\bigr)\bigr)
\]
is convex and its value at $\lambda=0$ is negative. Hence it can take
the value $0$ for a single value of $\lambda>0$ at most. When
$f(C)\equiv0$, the equation reduces to
\[
\e^{-\lambda}=\mathbb{E}_1(\exp(-\lambda M)).
\]
Recall that the Galton--Watson process is critical or sub-critical, so
$\mathbb{E}_1(M)\leq1$
according to Corollary 1 of \cite{Be2}. It is well known that this
ensures uniqueness of the solution to the preceding equation.
\end{pf}

Lemma \ref{L1} provides an implicit characterization of the law of the
partition into colonies through that of the pair of random variables
$(C,M)$. In turn, the latter can be conveniently described in terms of
a pair of random walks. This has its root in a key observation for
Galton--Watson processes that goes back to Harris~\cite{Harris}, and
will have an important role here for the analysis of asymptotic behaviors.
Specifically, consider
\[
S^\mathrm{h}_k=\xi^{\mathrm{h}}_1+\cdots+ \xi^{\mathrm{h}}_k - k
\quad \mbox{and}\quad  S^\mathrm{m}_k= \xi^{\mathrm{m}}_1+\cdots+\xi^{\mathrm
{m}}_k,\qquad
k\in\Z_+.
\]
Next define for every integer $j\geq0$ the first passage time
\[
\tau_j=\inf\{k \dvtx  S^\mathrm{h}_k=-j\}.
\]
We lift the following useful identity from Lemma 3 in \cite{Be2}.

\begin{lemma}\label{L2} The pair $(\tau_1, S^\mathrm{m}_{\tau_1})$
has the
same law as $(C,M)$.
\end{lemma}

We refer to Theorem 1(ii) in \cite{Be1} or to Proposition 1 in \cite
{Be2} for an explicit formula for this distribution which is obtained
by a combinatorial argument
and extends the well-known result of Dwass \cite{Dwass}
for the total population of Galton--Watson processes and passage times
of downward skip free random walks.
In this direction we also mention that the sequence of the atoms of the
partition into colonies has the same distribution under $\P_a$ as
%
%
\begin{equation}\label{E0}
(\tau_j-\tau_{j-1}\dvtx  1\leq j\leq\eta_a)  \qquad\mbox{with }
\eta_a=\inf\{j\dvtx  j-S^\mathrm{m}_{\tau_j}=a\}.
\end{equation}
This follows from Section 2 in \cite{Be1}; see in particular Lemma 4 there.
The interested reader may wish to provide an alternative proof of Lemma
\ref{L1} based on this representation and using the strong Markov
property for random walks in place of the extended branching property
for Galton--Watson processes.

\section{Random walks, L\'evy processes and passage times}\label{sec3}

As our main goal is to investigate limits of partitions into colonies,
Lemmas \ref{L1} and \ref{L2} suggest that we should study asymptotics
of first passage times in random walks,
which is the purpose of this section.
We first introduce the asymptotic regimes that we shall consider later
on, and develop some of their consequences for passage times of certain
random walks and L\'evy processes. Our starting point is a classical
result of convergence for rescaled Galton--Watson processes toward
continuous state branching processes (in short, CSBP) that we now
recall. We refer to the monograph \cite{DuLG} by Duquesne and Le Gall
for a complete account, including some terminology which will not be
defined here.

We consider for each integer $n$ a sequence $(\xi_{n,k}\dvtx  k\in\N)$ of
i.i.d. copies of some $\Z_+$-valued random variable with mean at most $1$
which should be thought of as the number of children of a typical
individual in
the $n$th population model.
The basic assumption is that there exists a sequence $\alpha(n)$ with
$\lim_{n\to\infty} \alpha(n)/n=\infty$ and a process $(X_t, t\geq0)$
such that
%
%
\begin{equation}\label{E1}
\frac{1}{n}\bigl( \xi_{n,1}+\cdots+\xi_{n,[\alpha(n)t]}-[\alpha(n)t]
\bigr) \Longrightarrow X_t
\end{equation}
for some (and then all) $t>0$, where the notation $\Longrightarrow$
refers to convergence in distribution as $n\to\infty$.
More precisely, \eqref{E1} then can be reinforced to
weak convergence on the space of c\`adl\`ag processes
endowed with Skorohod's topology; see, for instance, Theorem 16.4 in
\cite
{Kallenberg}. Moreover, the limit $X=(X_t\dvtx  t\geq0)$ is necessarily a
L\'evy process which has no negative jumps and does not drift to
$+\infty$,
that is $\mathbb{E}(X_t)\in[-\infty,0]$
[this follows from the requirement that $\mathbb{E}(\xi_{n,k})\leq1$
for every~$n$].

The law of the L\'evy process $X$ is characterized by its Laplace
exponent $\psi\dvtx  \R_+\to\R_+$ which is defined by
\[
\mathbb{E}(\exp(-qX_t))=\exp(t\psi(q)),\qquad  q\geq0.
\]
We shall further assume that $X$ has infinite variation, or
equivalently that
\[
\lim_{q\to\infty} q^{-1}\psi(q) =\infty.
\]
Next, consider a sequence $(a(n)\dvtx  n\in\N)$ with $a(n)/n\to a$ for
some $a>0$,
and denote by $Z^{(n)}$ a Galton--Watson process started from $a(n)$
ancestors and with reproduction law given by the distribution of $\xi_{n,k}$.
Then we have
\[
\frac{1}{n}Z^{(n)}_{[\alpha(n)t/n]} \Longrightarrow Z_t,
\]
where $(Z_t\dvtx  t\geq0)$ is a CSBP started from $Z_0=a$ and with
branching mechanism~$\psi$; see, for example, Theorem 2.1.1 in \cite{DuLG}.

We now turn our attention to the spatial case where some children of a
parent may emigrate.
That is we consider an array $((\xi^{\mathrm{h}}_{n,k}, \xi^{\mathrm
{m}}_{n,k})\dvtx  k,n\in\N)$ of random variables with values in $\Z_+^2$,
where $\xi^{\mathrm{h}}_{n,k}$ should be thought of as the number of
homebody children and $\xi^{\mathrm{m}}_{n,k}$ as the number of migrant
children of the $k$th individual for the $n$th population model; in particular
$\xi_{n,k}=\xi^{\mathrm{h}}_{n,k}+ \xi^{\mathrm{m}}_{n,k}$.
We assume that for each fixed integer $n$, the sequence $((\xi
^{\mathrm{h}}_{n,k}, \xi^{\mathrm{m}}_{n,k})\dvtx  k\in\N)$ is i.i.d.,
and just as in
the preceding section, we construct a pair of random walks
\[
S^\mathrm{h}_{n,k}=\xi^{\mathrm{h}}_{n,1}+\cdots+ \xi^{\mathrm
{h}}_{n,k} - k
\quad \mbox{and}\quad  S^\mathrm{m}_{n,k}= \xi^{\mathrm{m}}_{n,1}+\cdots+\xi
^{\mathrm{m}}_{n,k},\qquad
k\in\Z_+.
\]
Observe that the random walk $S^\mathrm{m}_{n,\cdot}$ is nondecreasing
while $S^\mathrm{h}_{n,\cdot}$ is downward skip free, that is, its increments
belong to $\{-1,0,1,2, \ldots\}$.
We now reinforce \eqref{E1} by assuming that the L\'evy process $X$ can
be decomposed as a sum
\[
X_t=X^\mathrm{h}_t+X^\mathrm{m}_t,
\]
where $((X^\mathrm{h}_t, X^\mathrm{m}_t)\dvtx  t\geq0)$ is a bivariate L\'
evy process,
in such a way that for some (and then all) $t>0$
%
%
\begin{equation}\label{E2}
\frac{1}{n}\bigl( S^\mathrm{h}_{n,[\alpha(n)t]}, S^\mathrm{m}_{n,[\alpha(n)t]}
\bigr) \Longrightarrow(X^\mathrm{h}_t,X^\mathrm{m}_t).
\end{equation}
We point out that again \eqref{E2} is automatically reinforced to weak
convergence in the sense of Skorohod by an appeal to Theorem 16.4 in
\cite{Kallenberg}.

We stress that necessarily, the L\'evy process $X^\mathrm{h}$ has no
negative jumps, infinite variation, and does not drift to $+\infty$,
and that $X^\mathrm{m}$ must be a subordinator (i.e., an increasing L\'evy
process); the two may or not be correlated.
We denote the bivariate Laplace exponent by $\Psi$, that is,
\[
\mathbb{E}\bigl(\exp-(q X^\mathrm{h}_t+rX^\mathrm{m}_t)\bigr)=\exp(t\Psi
(q,r)),\qquad
q,r\geq0.
\]
In particular, there is the identity
\[
\psi(q)=\Psi(q,q),\qquad  q\geq0 ;
\]
note also that our assumptions force $\Psi(q,q)\geq0$ whereas $\Psi
(0,r)\leq0$.

Next, we consider the first passage process
\[
T_x=\inf\{t\geq0\dvtx  X^\mathrm{h}_t<-x\},\qquad  x\geq0,
\]
which is a subordinator whose Laplace exponent is given by the inverse function
of $\Psi(\cdot,0)$; see Theorem VII.1 in \cite{LP}. Using $T_{\cdot}$
as a time-substitution, we also introduce
the compound process
\[
Y_x=X^\mathrm{m}_{T_x},\qquad  x\geq0.
\]
The distribution of the pair $(T,Y)$ can be described as follows.

\begin{lemma}\label{L3} \textup{(i)} The process
\[
\bigl((T_x, Y_x) \dvtx  x\geq0\bigr)
\]
is a bivariate subordinator.
\begin{longlist}[(iii)]
\item[(ii)] Its Laplace exponent $\Phi\dvtx  \R_+^2\to\R_+$
defined by
\[
\mathbb{E}\bigl(\exp(-qT_x-rY_x)\bigr)=\exp(-x\Phi(q,r)),\qquad
q,r\geq0,
\]
is determined as the unique solution to the equation
\[
\Phi(\Psi(q,r),r)=q.
\]

\item[(iii)]
There exists a unique measure $\Lambda$ on $\R_+^2\backslash\{(0,0)\}$
with $\int(1\wedge(x_1+x_2))\Lambda(\mathrm{d}x_1\, \mathrm
{d}x_2)<\infty$
such that
\[
\Phi(q,r)=\int(1-\e^{-qx_1-r x_2})\Lambda(\mathrm{d}x_1\, \mathrm
{d}x_2).
\]
In other words, the bivariate subordinator $(T_x, Y_x)$
has no drift and L\'evy measure~$\Lambda$.

\item[(iv)]
Finally, we also have
\[
\mathbb{E}(Y_1)= \int x_2 \Lambda(\mathrm{d}x_1\, \mathrm{d}x_2) \leq
1.
\]
\end{longlist}
\end{lemma}

\begin{pf} The proof is essentially a variation of that of Theorem VII.1 in
\cite{LP}.
The passage times $T_x$ are stopping times in the natural filtration of
the bivariate L\'evy process $(X^\mathrm{h},X^\mathrm{m})$ which are a.s.
finite and such that $X^\mathrm{h}_{T_x}=-x$ (by the absence of negative
jumps for $X^\mathrm{h}$ and the fact that $X^\mathrm{h}$ does not
drift to
$+\infty$). The strong Markov property immediately implies that
$(T_x, Y_x)$ is has independent and stationary increments; further this
process has clearly c\`adl\`ag nondecreasing sample paths in each
coordinate. In other words, it is a bivariate subordinator.

The Laplace exponent $\Phi$ is then determined by an application of
Doob's sampling theorem to the martingale
\[
\exp\bigl(-(qX^\mathrm{h}_t+rX^\mathrm{m}_t)-t\Psi(q,r)\bigr),\qquad  t\geq0
\]
(recall that $X^\mathrm{h}_{T_x}=-x$ a.s.). Observe that our assumptions
ensure that for every $r\geq0$, the function
$\Psi(\cdot,r)$ is continuous and convex with $\Psi(0,r)\leq0$ and
$\Psi(\infty,r)
=\infty$, so the equation $\Phi(\Psi(q,r),r)=q$ determines $\Phi$
on $\R_+^2$.

It remains to check that both subordinators have no drift. We know from
Corollary VII.5 in \cite{LP} and the fact that $X^\mathrm{h}$ has
unbounded variation that
$\lim_{q\to\infty}\Psi(q,0)/q=\infty$. This implies that
$\lim_{q\to\infty}\Phi(q,0)/q=0$ and hence $T_{\cdot}$ has no
drift. On
the other hand, the L\'evy--It\^o decomposition enables us to express
the subordinator $X^\mathrm{m}$ as the sum of a linear drift and a
pure-jump process. The time-substitution by $T_x$ thus yields that
$Y_x$ can be expressed as the sum of two pure-jump processes, and hence
its drift coefficient must be zero.

The penultimate displayed identity of the statement is just the
celebrated L\'evy--Khintchine formula. Finally, the assumption that the
L\'evy process $X$ does not drift to $+\infty$
is equivalent to requiring that its first moment exists and is nonpositive,
$\mathbb{E}(X_t)\leq0$. It follows that $X_t=X^\mathrm
{h}_t+X^\mathrm{m}_t$ is a
super-martingale,
and since $T_x$ is a stopping time, we deduce from Doob's sampling theorem
that for every $x,t\geq0$,
\[
\mathbb{E}(X^\mathrm{m}_{t\wedge T_x}) \leq\mathbb{E}(-X^\mathrm
{h}_{t\wedge T_x}) \leq x,
\]
where the second inequality is due to the definition of $T_x$ and the
absence of negative jumps for $X^\mathrm{h}$. Then it suffices to let
$t\to
\infty$ to get by monotone convergence that $\mathbb{E}(Y_x)\leq x$,
which in
turn yields our last claim by an application of the L\'evy--It\^o
decomposition of the subordinator $Y$
and the first-moment formula for Poisson measures.
\end{pf}

Lemmas \ref{L1} and \ref{L2} suggest that the asymptotic behavior of
the distribution of the partition into colonies should be related to
that of the first passage times
of the downward skip free random walk $S^\mathrm{h}_{n,\cdot}$,
\[
\tau_{n,j}=\inf\{k \dvtx  S^\mathrm{h}_{n,k}=-j\}, \qquad j\in\N.
\]
In this direction, we point at the following limit theorem.

\begin{corollary} \label{C1} In the regime \eqref{E2}, we have
for every bounded continuous function
$g\dvtx  \R_+^2\to\R$ with $g(x_1,x_2)=O(x_1+x_2)$ as $x_1+x_2\to0$ that
\[
\lim_{n\to\infty} n \mathbb{E}(g(\alpha(n)^{-1} \tau_{n,1},
n^{-1} S^\mathrm{m}_{n,\tau_{n,1}})) =\int_{\R
_+^2}g(x_1,x_2)\Lambda(\mathrm{d}x_1\,
\mathrm{d}x_2),
\]
where the L\'evy measure
$\Lambda$ has been defined in Lemma \ref{L3}.
\end{corollary}

\begin{pf} It follows from the assumptions \eqref{E2} and routine arguments
[recall that~$ S^\mathrm{h}_{n,\cdot}$ is downwards skip free and that
\eqref{E2} can be reinforced to weak convergence of c\`adl\`ag
processes]
that for an arbitrary $x>0$
%
%
\begin{equation}\label{E3'}
\bigl(\alpha(n)^{-1}\tau_{n,[nx]},n^{-1}
S^\mathrm{m}_{n,\tau_{n,[nx]}}\bigr) \Longrightarrow (T_x, Y_x).
\end{equation}
On the other hand, one readily deduces from the strong Markov property
for random walks that for each fixed $n$,
\[
(\tau_{n,k},
S^\mathrm{m}_{n,\tau_{n,k}}),\qquad   k\geq0,
\]
is a random walk with nondecreasing coordinates.
We complete the proof by taking $x=1$ in \eqref{E3'}, and appealing to
(i) in Corollary 15.16 in \cite{Kallenberg} and Lemma \ref{L3}.
\end{pf}

\section{A family of random point measures}\label{sec4}

In this section, we introduce and develop some properties of a class of
random point measures which will arise later on as limits for
partitions into colonies.
The idea stems from the representation \eqref{E0} of the sequence of
the atoms of the partition into colonies. Indeed, as by the strong
Markov property, the increments
$\tau_j-\tau_{j-1}$ of the first passage time process in a downward
skip free random walk
are i.i.d., the combination of \eqref{E0} and the law of rare events
suggest that if a limiting partition exists, then it should be
described in terms of a Poisson random measure restricted to a random
domain with a boundary given by a first passage time.

Our basic analytic datum is some sigma-finite measure on $\R_+^2$ with
no mass at $(0,0)$ that will be denoted by $\Lambda$. Although in
subsequent sections $\Lambda$ will be chosen to be the L\'evy measure
that arises in Lemma \ref{L3}, this specification is not required in
the present section
(of course, the notation introduced here is coherent with that of
Section \ref{sec3}). We write $\Lambda^{(1)}$ and $\Lambda^{(2)}$ for the restrictions
to $(0,\infty)$ of the two marginals of $\Lambda$ and assume that
$\Lambda^{(1)}$ is sigma-finite and
%
%
\begin{equation}\label{E3}
\int_{(0,\infty)} x\Lambda^{(2)}( \mathrm{d} x)\leq1.
\end{equation}

We consider a Poisson measure ${\mathcal N}$ on $(0,\infty)\times\R_+^2$
with intensity measure $\mathrm{d}t\otimes\Lambda(\mathrm{d}x_1\,
\mathrm{d}
x_2)$, and denote by $(t,\Delta_t)=(t,\Delta^{(1)}_t,\Delta
^{(2)}_t)$ a generic
atom of ${\mathcal N}$. Following the classical construction of L\'evy
and It\^o, we introduce the subordinator
\[
Y_t
=\int_{(0,t]\times\R_+^2}x_2 {\mathcal N}(\mathrm{d}s\, \mathrm
{d}x_1 \,\mathrm{d}x_2)
=\sum_{0<s\leq t}\Delta^{(2)}_s,\qquad   t\geq0.
\]
Because $Y$ has no drift and its L\'evy measure $\Lambda^{(2)}$ fulfills
\eqref{E3}, we have $\mathbb{E}(Y_1)\leq1$.
In particular, the first passage times
\[
\sigma_y=\inf\{t\geq0\dvtx  t-Y_t=y\},\qquad  y\geq0,
\]
are finite a.s.

Next, for every $t\geq0$, we consider the point measure on $(0,\infty)$
\[
{\mathcal N}^{(1)}_t=\sum_{0< s\leq t}\delta_{\Delta^{(1)}_s}.
\]
Note that ${\mathcal N}^{(1)}_t$ is a Poisson random measure with intensity
$t\Lambda^{(1)}$, and from the superposition property of Poisson
measure, that
the measure-valued process
$({\mathcal N}^{(1)}_t\dvtx  t\geq0)$ has independent and stationary increments.
The random point measures we are interested in are defined
by time-substitution through the passage times $\sigma_y$
%
%
\begin{equation}\label{E4}
{\mathcal M}_y={\mathcal N}^{(1)}_{\sigma_y}=\sum_{0<t\leq\sigma
_y}\delta
_{\Delta^{(1)}_t}.
\end{equation}
In words, ${\mathcal M}_y $ is the image of the restriction of
${\mathcal N}$ to
the random set $(0,\sigma_y]\times\R_+^2$ by the projection
$(s,x_1,x_2)\to x_1$.

In the special case when the intensity measure $\Lambda$ is carried by
the diagonal
$\{(x,x)\dvtx  x>0\}$, that is, when $\Lambda(\mathrm{d}x_1\, \mathrm
{d}x_2)=\delta
_{x_1}(\mathrm{d}x_2)
\Lambda^{(1)}(\mathrm{d}x_1)$, we have $\Delta^{(1)}_s=\Delta
^{(2)}_s$ a.s. and the
random point measure ${\mathcal M}_y $ coincides with the empirical
measure of the sizes of the
jumps performed by the subordinator $Y$ during the time-interval
$(0,\sigma_y]$. This case has also a natural interpretation in terms of
continuous state branching processes in discrete time (see \cite{Jirina}).
More precisely,
the well-known correspondence between CSBP in discrete time
and subordinators enables us to think of ${\mathcal M}_y $ as the
empirical measure of the sizes of siblings in a CSBP in discrete time
with reproduction intensity $\Lambda^{(1)}$ and started from an initial
population of size~$y$.

We now observe that the property of independence and stationarity of
the increments for the process of point measures
$({\mathcal N}^{(1)}_t\dvtx  t\geq0)$ is preserved after the time-substitution
by $\sigma_y$. This claim is essentially a variation of the well-known
fact that
the first passage process of a real-valued L\'evy process with no
negative jumps is a subordinator; see, for example, Theorem VII.1 in \cite{LP}.

\begin{lemma}\label{L4}
The measure-valued process $({\mathcal M}_y\dvtx  y\geq0)$ has independent
and stationary increments.
\end{lemma}

\begin{pf} Assume for a while that $\int_{(0,\infty)}(1\wedge x)\Lambda^{(1)}
(\mathrm{d}x)<\infty$,
which enables us to construct
\[
T_t=\bigl\langle{\mathcal N}^{(1)}_t, \mathrm{Id}\bigr\rangle
=\int_{(0,t]\times\R_+^2}x_1 {\mathcal N}(\mathrm{d}s\, \mathrm
{d}x_1\, \mathrm{d}x_2)
=\sum_{0<s\leq t}\Delta^{(1)}_s, \qquad t\geq0.
\]
Plainly $(T,Y)$ is a pure-jump L\'evy process, more precisely it is a
bivariate subordinator with no drift. Further, the Poisson measure
${\mathcal N}$ can be recovered from
the jump process of $(T,Y)$.
The strong Markov property for L\'evy processes shows that
for every $y>0$, the shifted process
\[
(T',Y')_t=(T,Y)_{\sigma_y+t}-(T,Y)_{\sigma_y}, \qquad t\geq0,
\]
is independent of $((T,Y)_t\dvtx  t\leq\sigma_y)$
and has the same law as $(T,Y)$.
As $\sigma_{y+y'}-\sigma_y$ coincides with the first passage time of
the process $t\to t-Y'_t$ at level $y'$, this establishes our claim.

Finally, the assumption that
$\int_{(0,\infty)}(1\wedge x)\Lambda^{(1)}(\mathrm{d}x)<\infty$
can be removed
by considering
the image of ${\mathcal N}$ by a mapping $(t,x_1,x_2)\to(t, \phi(x_1),
x_2)$ for some
appropriate bijective map $\phi$ (recall that the measure $\Lambda^{(1)}$
is sigma-finite).
\end{pf}

We next point at an interesting connection between the distributions
of ${\mathcal N}^{(1)}_t$ and~${\mathcal M}_y$ which is an avatar of the
classical ballot theorem (compare with Corollary~VII.3 in \cite{LP}).

\begin{proposition}\label{P1} There is the identity
\[
\P({\mathcal M}_y \in A, \sigma_y\in\mathrm{d}t)\,\mathrm{d}y
=\frac{y}{t}\P\bigl({\mathcal N}^{(1)}_t\in A, t-Y_t\in\mathrm
{d}y\bigr)\,\mathrm{d}t
,\qquad t>y>0,
\]
where $A$ denotes an arbitrary measurable subset of point measures on
$(0,\infty)$.
\end{proposition}
\begin{pf} Introduce the random set
\[
{\mathcal R}=\Bigl\{t\dvtx  t-Y_t=\max_{0\leq s \leq t} (s-Y_s)\Bigr\}.
\]
The cyclic exchangeability property of the point measure ${\mathcal N}$
enables us
to use a variation of the well-known combinatorial argument for the
ballot theorem (see~\cite{Takacs}) and get
\[
t\mathbb{E}\bigl(g(t-Y_t),  {\mathcal N}^{(1)}_t\in A\mbox{ and }t\in
{\mathcal R}\bigr)
= \mathbb{E}\bigl(g(t-Y_t) (t-Y_t)^+
, {\mathcal N}^{(1)}_t\in A\bigr),
\]
where $g\dvtx \R\to[0,\infty)$ stands for a generic measurable function.
This easily yields the claim.
\end{pf}

\begin{remark*} In the special when the intensity measure
$\Lambda$ is carried by the diagonal, we have $Y_t=\langle{\mathcal
N}^{(1)}_t, \mathrm{Id}\rangle$ a.s.,
and Proposition \ref{P1} shows that the distribution of ${\mathcal
M}_y$ is essentially a mixture of laws of \textit{Poisson--Kingman
partitions} as defined by Pitman \cite{PiPK}. More precisely, suppose
for simplicity that for every $t>0$, the infinitely divisible variable
$Y_t$ has an absolutely continuous law with a continuous density, say
$\rho_t(\cdot)$.
It then follows from Proposition \ref{P1} that
\[
\P({\mathcal M}_ y\in\cdot)
= y\int_{(y,\infty)}\frac{1}{t}\P\bigl({\mathcal N}^{(1)}_t\in\cdot
\mbox{ }\big|\mbox{ }
\bigl\langle{\mathcal N}^{(1)}_t, \mathrm{Id}\bigr\rangle=t-y\bigr) \rho_t(t-y)
\,\mathrm{d}t,
\]
where ${\mathcal N}^{(1)}_t$ is a Poisson random measure with intensity
$t\Lambda^{(1)}$.
Up-to a normalization, the conditional Poisson measures
$\P({\mathcal N}^{(1)}_t\in\cdot\mid\langle{\mathcal N}^{(1)}_t,
\mathrm{Id}\rangle=a) $
which appear in the integral above belong to the family of
Poisson--Kingman partitions studied in depth by Pitman \cite{PiPK}.
\end{remark*}

Next, for every Borel function $f\dvtx (0,\infty)\to\R_+$ with compact
support, we define
the cumulant $\kappa(f)>0$ by
\[
\mathbb{E}(\exp-\langle{\mathcal M}_1,f\rangle)=\exp(-\kappa(f)),
\]
with the usual notation
\[
\langle{\mathcal M}_1,f\rangle=\sum_{0<t\leq\sigma_1}f\bigl(\Delta
^{(1)}_t\bigr).
\]
Observe from Lemma \ref{L4} that for an arbitrary $y>0$ we have more generally
%
%
\begin{equation}\label{E5}
\mathbb{E}(\exp-\langle{\mathcal M}_y,f\rangle)=\exp(-y\kappa(f)).
\end{equation}
It is well known that the cumulant $\kappa$ determines the law of
${\mathcal M}_1$, in the sense that any random measure on $(0,\infty)$
having the same cumulant as ${\mathcal M}_1$ is distributed as~${\mathcal M}_1$; see for instance Lemma 12.1 in \cite
{Kallenberg}. We may now state the following basic result which
provides the characteristic equation solved by the cumulant:

\begin{theorem} \label{T1}
For every Borel function $f\dvtx [0,\infty)\to\R_+$ with compact support in
$(0,\infty)$, the equation
\[
\lambda=\int_{\R_+^2}\bigl(1-\exp\bigl(-f(x_1)-\lambda x_2\bigr)\bigr)\Lambda
(\mathrm{d}x_1\, \mathrm{d} x_2)
\]
has a unique solution in $[0,\infty)$ which is given by $\lambda
=\kappa(f)$.
\end{theorem}

\begin{pf}
For any random time $R\geq1$, we see from elementary properties of
Poisson random measures that ${\mathcal N}^{(1)}_R={\mathcal
N}^{(1)}_1+\tilde
{\mathcal N}^{(1)}_{R-1}$
where $(\tilde{\mathcal N}_t\dvtx  t\geq0)$ is a process of point
measures which is independent of
the restriction of ${\mathcal N}$ to $[0,1]\times\R_+^2$ and has the
same distribution as $({\mathcal N}_t\dvtx  t\geq0)$. We then note that the
first passage time~$\sigma_1$ is bounded from below by $1$, and more
precisely there is the identity
\[
\sigma_1=1+\tilde\sigma(Y_1),
\]
where
\[
\tilde\sigma(y)=\inf\{t\geq0\dvtx  t-\tilde Y_t=y\}\quad  \mbox{and}\quad
\tilde Y_t=\int_{[0,t]\times\R_+^2}x_2 \tilde{\mathcal N}(\mathrm{d}s\,
\mathrm{d}x_1\, \mathrm{d}x_2).
\]
Applying the preceding observation, we thus have
%
%
\begin{equation}\label{E6}
{\mathcal M}_1={\mathcal N}^{(1)}_{\sigma_1}={\mathcal
N}^{(1)}_1+\tilde
{\mathcal N}^{(1)}_{\tilde\sigma(Y_1)},
\end{equation}
so we can deduce from \eqref{E5} that
\[
\exp(-\kappa(f)) =
\mathbb{E}\bigl(\exp-\bigl(\bigl\langle{\mathcal N}^{(1)}_{1},f\bigr\rangle+ Y_1 \kappa(f)\bigr)
\bigr).
\]
From the very definitions of ${\mathcal N}^{(1)}$ and $Y_1$, we can
rewrite the preceding identity as
\begin{eqnarray*}
\exp(-\kappa(f))&=&\mathbb{E}\biggl(\exp\biggl(-\sum_{0<t\leq1}\bigl( f\bigl(\Delta^{(1)}
_t\bigr)+\kappa(f)\Delta^{(2)}_t\bigr)\biggr)
\biggr) \\
&=&\exp\biggl(-\int_{\R_+^2}\bigl(1-\exp\bigl(-f(x_1)-\kappa(f) x_2\bigr)
\bigr)\Lambda(\mathrm{d}x_1 \,\mathrm{d} x_2)\biggr),
\end{eqnarray*}
where the last line is Campbell's identity. Thus $\kappa(f)$ solves
\eqref{E4}.

Uniqueness is now easy. Indeed the map
\[
F\dvtx \lambda\to\lambda-\int_{\R_+^2}\bigl(1-\exp\bigl(-f(x_1)-\lambda
x_2\bigr)\bigr)\Lambda(\mathrm{d}x_1\, \mathrm{d} x_2)
\]
has derivative
\[
F'(\lambda)=1-\int_{\R_+^2} x_2\exp\bigl(-f(x_1)-\lambda x_2\bigr)\Lambda
(\mathrm{d}x_1 \,\mathrm{d} x_2)
\]
which is positive due to \eqref{E3}.
\end{pf}

We now conclude this section by discussing a simple example. Suppose that
$\Lambda$ has support on the axes, that is,
\[
\Lambda(\mathrm{d}x_1\, \mathrm{d} x_2)
=\Lambda^{(1)}(\mathrm{d}x_1)\delta_0 (\mathrm{d} x_2) + \delta
_0(\mathrm{d}x_1)\Lambda^{(2)}( \mathrm{d} x_2).
\]
Then the equation in Theorem \ref{T1} can be rewritten as
\[
\int_{(0,\infty)}\bigl(1-\e^{-f(x_1)}\bigr)\Lambda^{(1)}(\mathrm{d}x_1)
=\lambda-\int_{(0,\infty)}(1-\e^{-\lambda x_2})\Lambda^{(2)}
(\mathrm{d} x_2).
\]
On the other hand, our assumption implies that the subordinator $Y$ and
the process of point measures ${\mathcal N}^{(1)}$ are independent, and
${\mathcal M}_1={\mathcal N}^{(1)}_{\sigma_1}$ is thus a mixed Poisson
measure with
intensity $t \Lambda^{(1)}$ and mixing law $\P(\sigma_1\in\mathrm
{d}t)$. In
particular we have
\[
\mathbb{E}(\exp-\langle{\mathcal M}_1,f\rangle)=
\int_{(0,\infty)}\exp\biggl(-t\int_{(0,\infty)}\bigl(1-\e^{-f(x_1)}\bigr)\Lambda^{(1)}
(\mathrm{d}x_1)\biggr)
\P(\sigma_1\in\mathrm{d}t).
\]
Now recall from Theorem VII.1 in \cite{LP} that the Laplace transform
of the first passage time $\sigma_1$ of the L\'evy process with no
positive jumps
$t-Y_t$ is given by
\[
\mathbb{E}(\e^{-q\sigma_1})=\exp(-\varphi(q)),\qquad   q\geq0,
\]
where the cumulant $\varphi$ is the unique solution to
\[
q=\varphi(q)-\int_{(0,\infty)}\bigl(1-\e^{-\varphi(q)x_2}\bigr)\Lambda
^{(2)}(\mathrm{d}x_2).
\]
We conclude that
\[
\kappa(f) = \varphi\biggl(\int_{(0,\infty)}\bigl(1-\e^{-f(x_1)}\bigr)\Lambda^{(1)}
(\mathrm{d}x_1)\biggr),
\]
which is thus in agreement with Theorem \ref{T1}.

\section{Limit laws for partitions into colonies}\label{sec5}

In this section, we state and prove the main limit theorem for
distributions of partitions into colonies. We consider for each fixed
integer $n$ a Galton--Watson process with emigration started from~$a(n)$ ancestors that all occupy different sites, such that the number
of homebody children and the number of migrant children $(\xi^{\mathrm
{h}}_{n,k},\xi^{\mathrm{m}}_{n,k})$ of the $k$th individual is given
by an
i.i.d. sequence.
We write ${\mathcal P}_n$ for the partition into colonies induced by
this model.

We also consider a bivariate L\'evy process $(X^\mathrm{h},X^\mathrm{m})$
such that $X^\mathrm{h}$ has no negative jumps and infinite variation,
$X^\mathrm{m}$ is a subordinator, and the sum $X=X^\mathrm
{h}+X^\mathrm{m}$ does
not drift to $+\infty$. We write $\Lambda$ for the L\'evy measure that
arises in Lemma \ref{L3},
and then $({\mathcal M}_y\dvtx  y\geq0)$ for the process of random point
measures which has been studied in Section \ref{sec4} for this specific choice
of $\Lambda$.

We are now able to state the main result of this work.

\begin{theorem}\label{T2} Write $\tilde{\mathcal P}_n$ for the image of
the partition into colonies ${\mathcal P}_n$ by the rescaling $x\to
x/\alpha(n)$, namely,
\[
\langle\tilde{\mathcal P}_n, f\rangle= \langle{\mathcal P}_n,
\tilde f_n\rangle,
\]
where $\tilde f_n(x) = f(x/\alpha(n))$.
Assume that the number of ancestors $a(n)$ fulfills $a(n)\sim a n$ for
some $a>0$.
Then in the regime \eqref{E2}, $\tilde{\mathcal P}_n$ converges weakly
on the space of sigma-finite measures on $(0,\infty)$ as $n\to\infty$
toward ${\mathcal M}_a$.
\end{theorem}

The material developed so far suggests that the proof of Theorem \ref
{T2} should consist of two steps, namely first a tightness property for
the rescaled partitions into colonies,
and then uniqueness of the limit of a subsequence that shall be derived
by the analysis of cumulants. This is indeed the route that we will follow.

\begin{lemma} \label{L5} In the regime \eqref{E2}, the sequence of the
distributions of the variables $\langle\tilde{\mathcal P}_n, \mathrm
{Id}\rangle$, for $n\in\N$,
is tight on the space of sigma-finite measures on $(0,\infty)$.

\end{lemma}
\begin{pf}
Indeed, recall that
\[
\langle\tilde{\mathcal P}_n, \mathrm{Id}\rangle= \alpha(n)^{-1}
\langle{\mathcal P}_n, \mathrm{Id}\rangle=\alpha(n)^{-1}\sum
_{j=1}^{\gamma
^{(n)}} C_j^{(n)} = \zeta_n/\alpha(n)
\]
is simply the size $\zeta_n$ of the total population generated by the
Galton--Watson process $Z^{(n)}$ renormalized by the factor $1/\alpha
(n)$. It is well known
that in the regime~\eqref{E2}, this quantity converges in distribution
as $n\to\infty$ toward the size of the total population of the CSBP
$Z$, that is, equivalently,
the first passage time of the L\'evy process $X$ at level $-a$.
Hence, the sequence in the statement is tight.
\end{pf}

For the second step of the proof of Theorem \ref{T2}, we write
$\tilde K_n$ for the cumulant of the rescaled random measure $\tilde
{\mathcal P}_n$, and fix a Borel function $f\dvtx (0,\infty)\to\R_+$ with
compact support. For the sake of simplicity, we will suppose in the
sequel that $a=1$, that is, that $a(n)\sim n$, which induces no loss of
generality thanks to the branching property and \eqref{E5}.

\begin{lemma} \label{L6} Let $f\dvtx [0,\infty)\to\R_+$ be an arbitrary
continuous function with compact support in $(0,\infty)$.
Any limit point $\lambda$ of the sequence $(\tilde K_n(f)\dvtx n\in\N)$
fulfills
\[
\lambda=\int_{\R_+^2}\bigl(1-\exp\bigl(-f(x_1)-\lambda x_2\bigr)\bigr)\Lambda
(\mathrm{d}x_1 \,\mathrm{d} x_2).
\]
\end{lemma}

\begin{pf} We work with an increasing subsequence of integers $n$ such that
$\tilde K_n(f)\to\lambda$.
Tracing back the definitions, we get
\begin{eqnarray*}
\exp(-\tilde K_n(f))&=&\exp\bigl(-a(n)K_n(\tilde f_n)\bigr)\\
&=&\bigl( 1-\mathbb{E}^{(n)}_1\bigl(1-\exp\bigl(-\tilde f_n(C)- K_n(\tilde f_n) M\bigr)\bigr)
\bigr)^{a(n)},
\end{eqnarray*}
where we used Lemma \ref{L1} for the last equality, and the notation
$\mathbb{E}
^{(n)}$ refers to the mathematical expectation corresponding to the
$n$th population model. Taking logarithms, we arrive at
\[
a(n) \mathbb{E}^{(n)}_1\bigl(1-\exp\bigl(-\tilde f_n(C)- K_n(\tilde f_n) M\bigr)\bigr)
\to\lambda.
\]

Then recall from Lemma \ref{L2} that
\begin{eqnarray*}
& &\mathbb{E}^{(n)}_1\bigl(1-\exp\bigl(-\tilde f_n(C)- K_n(\tilde f_n) M\bigr)\bigr)\\
&&\qquad=\mathbb{E}^{(n)}_1\biggl(1-\exp\biggl(-f\bigl(C/\alpha(n)\bigr)- \frac{n}{a(n)} \tilde
K_n(f) \frac{M}{n}\biggr)\biggr)\\
&&\qquad=\mathbb{E}\biggl(1-\exp\biggl(-f(\alpha(n)^{-1}\tau_{n,1})- \frac{n}{a(n)}
\tilde
K_n(f) n^{-1} S^\mathrm{m}_{n,\tau_{n,1}}\biggr)\biggr).
\end{eqnarray*}
Recall also that ${n}/{a(n)}\to1$ and $\tilde K_n(f)\to\lambda$. We
now see that for every $\varepsilon>0$
\[
\limsup n \mathbb{E}\bigl(1-\exp\bigl(-f(\alpha(n)^{-1}\tau_{n,1})- (\lambda
+\varepsilon) n^{-1} S^\mathrm{m}_{n,\tau_{n,1}}\bigr)\bigr)\leq\lambda,
\]
so applying Corollary \ref{C1} with $g(x_1,x_2)=1-\exp
(-f(x_1)-(\lambda
+ \varepsilon) x_2)$, we get
\[
\int_{\R_+^2}\bigl(1-\exp\bigl(-f(x_1)-(\lambda+\varepsilon) x_2\bigr)
\bigr)\Lambda(\mathrm{d}x_1\, \mathrm{d} x_2)\leq\lambda.
\]
By a similar argument, we also obtain
\[
\lambda\leq\int_{\R_+^2}\bigl(1-\exp\bigl(-f(x_1)-(\lambda-\varepsilon)
x_2\bigr)\bigr)\Lambda(\mathrm{d}x_1\, \mathrm{d} x_2).
\]
We derive the equation of the statement letting $\varepsilon$ tend to
$0$.
\end{pf}

The proof of Theorem \ref{T2} should now be plain. It follows from
Lemma \ref{L5}
that the sequence of the laws of rescaled partitions into colonies
$\tilde{\mathcal P}_n$ is tight in the space of sigma-finite measures
on $(0,\infty)$. We then deduce from Prohorov's lemma (see, e.g., Lemma
16.15 in \cite{Kallenberg}) that the sequence of the
distributions of the random measures $\tilde{\mathcal P}_n$ on
$(0,\infty)$ is relatively compact. If $\tilde{\mathcal P}$ has the law
of the limit of some sub-sequence, then we deduce from Lemma \ref{L6}
that for an arbitrary continuous function $f\dvtx [0,\infty)\to\R_+$ with
compact support in $(0,\infty)$,
the cumulant
\[
\tilde K(f)=-\ln\mathbb{E}(\exp-\langle\tilde{\mathcal P}
,f\rangle)
\]
solves
\[
\tilde K(f) =\int_{\R_+^2}\bigl(1-\exp\bigl(-f(x_1)-\tilde K(f) x_2\bigr)
\bigr)\Lambda(\mathrm{d}x_1\, \mathrm{d} x_2).
\]
We conclude from Theorem \ref{T1} that $\tilde K(f)=\kappa(f)$, and
thus $\tilde{\mathcal P}$ has the same distribution as ${\mathcal M}_1$.

\begin{remark*} It may be interesting to point at a different route for
establishing Theorem \ref{T2}, which uses the representation \eqref{E0}
of the partition into colonies.
Recall the notation there and the convergence in distribution \eqref
{E3'}. Invoking Theorem 16.14 in \cite{Kallenberg}, it is easy to check
that the latter can be reinforced into weak convergence of c\`adl\`ag
processes in the sense of Skorohod.
One can then deduce from a time-substitution that
\[
\bigl( \alpha(n)^{-1}\tau_{n,[nx]}, n^{-1}S^\mathrm{m}_{n,\tau
_{n,[nx]}}\bigr)
\Longrightarrow (T_x,Y_x),
\]
where again the convergence holds in the sense of Skorohod.
Loosely speaking, this entails the weak convergence of the increments of
the random walk $\tau_{n,\cdot}$ rescaled by a factor $1/\alpha(n)$ to
the jump-process of the subordinator $T$.
We know from the L\'evy--It\^o decomposition that the latter can be
described as a Poisson random measure whose intensity is expressed in
terms of the L\'evy measure
of $T$. It remains to recall that in this setting, the number $\gamma
_n$ of colonies fulfills
\[
\gamma_n=\min\{k\dvtx  S^\mathrm{m}_{\tau_{n,k}}-k=-a(n)\}
\]
and to check that
\[
n^{-1}\gamma_n \Longrightarrow\sigma_a=\inf\{t\geq0\dvtx  t-Y_t=a\}.
\]
Some technical details needed to justify rigorously this approach may
be tedious; they are circumvented here by the appeal to the
characterization of the cumulant of the random measure ${\mathcal M}_a$
in Theorem \ref{T1} and the simple argument for tightness in Lemma
\ref{L5}.
\end{remark*}

\section{Examples}\label{sec6}

In this section, we shall illustrate our main results for partitions
into colonies by discussing some natural examples. Their common feature
is that the distribution of the total number of children (homebodies
and migrants) of a typical individual is fixed, that is
the Galton--Watson process for which spatial locations of individuals
are discarded has a fixed reproduction law. The differences in the
models thus only appear through the repartition between homebody and
migrant children. One could, of course, deal with much more general
examples, however the present ones already exhibit a rich variety of
asymptotic behaviors.

Recall \eqref{E1}. Throughout this section, we consider an
integer-valued random variable $\xi$ with unit mean,
which belongs to the domain of attraction of a (completely asymmetric)
stable variable
with index $\beta\in(1,2]$.
That is, there is a sequence $(\alpha(n)\dvtx  n\in\N)$ which varies
regularly with index $\beta$ such that
\[
n^{-1}\bigl(\xi_1+\cdots+ \xi_{\alpha(n)}-\alpha(n)
\bigr)\Longrightarrow
X_1,
\]
where $(\xi_i\dvtx  i\in\N)$ is a sequence of i.i.d. copies of $\xi$ and now
$(X_t\dvtx  t\geq0)$ a $\operatorname{stable}(\beta)$ L\'evy process with no negative
jumps. In other words, there is some $b>0$ such that
\[
\mathbb{E}(\exp(-q X_t))=\exp(t b q^{\beta}),\qquad   q\geq0,
\]
i.e. $\psi(q)=bq^{\beta}$. The (continuous version of the) density of
the variable $X_1$ will be denoted by
$\rho$,
\[
\P(X_1\in\mathrm{d } x)= \rho(x)\, \mathrm{d } x,
\]
so that, by scaling,
\[
\P(X_t\in\mathrm{d } x)= t^{-1/\beta}\rho(t^{-1/\beta}x)\, \mathrm
{d } x
\]
for every $t>0$.

\subsection{Allelic partitions for rare neutral mutations}\label{sec6.1}
We first deal with the classical model corresponding to neutral mutations.
That is for each fixed integer $n$, the total number of children of the
$k$th individual
is decomposed as $\xi_k=\xi^\mathrm{h}_{n,k}+\xi^\mathrm{m}_{n,k}$ where
conditionally on
$\xi_k=\ell$, the variable $\xi^\mathrm{m}_{n,k}$ has the binomial
distribution with parameter
$(\ell,p(n))$ for some $p(n)\in(0,1)$. In other words, we assume that
each child chooses to become a migrant with probability $p(n)$,
independently of the other individuals.

If we now suppose that
\[
p(n) \sim c n/\alpha(n)
\]
for some constant $c>0$, so that the mutation rate is small when $n$ is large,
then \eqref{E2} clearly holds with
$X^\mathrm{h}_t=X_t-ct$ and $X^\mathrm{m}_t=ct$, and thus
\[
\Psi(q,r)=bq^{\beta}+cq -cr, \qquad q,r\geq0.
\]
We now see from Lemma \ref{L3}(ii) that the bivariate subordinator
$(T_x,Y_x)$ has
Laplace exponent $\Phi(q,r)=\varphi(q+cr)$ where $\varphi(q)=z$ is
given by the nonnegative solution to
\[
bz^{\beta}+cz =q.
\]
We stress that $Y_x=cT_x$ a.s., and that the subordinator $(T_x\dvtx  x\geq
0)$ has Laplace exponent $\varphi$.

An easy consequence of a version of the Ballot theorem (more precisely,
cf. Corollary VII.3 in \cite{LP})
is that the L\'evy measure $\Lambda^{(1)}$ of the subordinator of the first
passage time of $X^\mathrm{h}$ can be expressed in terms of the
density of
$X^\mathrm{h}$ at $0$. More precisely,
one gets
\begin{eqnarray*}
\Lambda^{(1)}( \mathrm{d } t) &=& \lim_{x\to0+} x^{-1}\P(T_x\in
\mathrm{d}t)\\
&=& \lim_{x\to0+} \frac{1}{t} \frac{\P(-X^\mathrm{h}_t\in\mathrm
{d }
x)}{\mathrm{d}x} \,\mathrm{d}t\\
&=& t^{-1-1/\beta}\rho(ct^{1-1/\beta}) \,\mathrm{d } t,
\end{eqnarray*}
where the last equality follows from the fact that
\[
\P(X^\mathrm{h}_t\in\mathrm{d } x)
=\P(X_t\in ct +\mathrm{d } x)= t^{-1/\beta}\rho\bigl(t^{-1/\beta
}(x+ct)\bigr) \,\mathrm{d
} x.
\]
Recall that $Y=cT$; it follows that the L\'evy measure of the bivariate
subordinator
$(T,Y)$ that determines the law of the limiting partition ${\mathcal
M}_a$ is then given by
\[
\Lambda(\mathrm{d}x_1\, \mathrm{d}x_2) =x_1^{-1-1/\beta}\rho
(cx_1^{1-1/\beta})
\delta_{cx_1}(\mathrm{d}x_2)\, \mathrm{d}x_1.
\]

On the other hand, recall again from Corollary VII.3 in \cite{LP} that
\[
\frac{\P(T_x\in\mathrm{d}t)}{\mathrm{d}t}
= \frac{x}{t} \frac{\P(-X^\mathrm{h}_t\in\mathrm{d } x)}{\mathrm{d}x}=
xt^{-1-1/\beta}\rho\bigl(t^{-1/\beta}(ct-x)\bigr), \qquad t,x>0.
\]
We can combine this identity with the argument in the remark following
Proposition \ref{P1} to express the distribution of ${\mathcal M}_a$ as
a mixture of laws of Poisson measures conditioned on their first
moments (i.e., Poisson--Kingman partitions; see \cite{PiPK}). More
precisely, we have
\[
\P({\mathcal M}_ a\in\cdot)
= a\int_{(a,\infty)}\frac{1}{t}\P\bigl({\mathcal N}^{(1)}_t\in\cdot
\mid
t-Y_t=a\bigr) \frac{\P(t-Y_t\in\mathrm{d}a)}{\mathrm{d} a} \,\mathrm{d}t,
\]
where ${\mathcal N}^{(1)}_t$ is a Poisson random measure with intensity
$t\Lambda^{(1)}$.
We now get, using the identity $Y=cT$, that
%
%
\begin{eqnarray} \label{E9}
 \P({\mathcal M}_ a\in\cdot)&= &a\int_{(a,\infty)}\P\bigl({\mathcal N}^{(1)}_t\in\cdot\mid
 T_t=(t-a)/c\bigr)
 \nonumber
 \\[-8pt]
 \\[-8pt]
 \nonumber
&&\quad \qquad{}\times\rho\bigl(t^{-1/\beta}\bigl(ct-(t-a)/c\bigr)\bigr) \frac{t-a}{c^2} t^{-2-1/\beta}
\,\mathrm{d}t.
\end{eqnarray}

In the important case $\beta=2$, which occurs whenever the reproduction
law has finite variance, $\rho$ is simply the Gaussian density and
\[
\Lambda^{(1)}( \mathrm{d } t)= \frac{1}{2\sqrt{\pi t^3 b}}\exp
\biggl(-\frac{c^2
t}{4b}\biggr)\, \mathrm{d } t.
\]
That is $\Lambda^{(1)}$ is the L\'evy measure of an inverse Gaussian
subordinator, which is merely an exponential transform of the
$\operatorname{stable}(1/2)$ L\'evy measure. Recall that the exponential transform plays
no role for the distribution of Poisson measures conditioned on their
first moments, which thus reduces the description \eqref{E9} of the law
of~${\mathcal M}_ a$ to the more usual $\operatorname{stable}(1/2)$ L\'evy measure. This
situation has been investigated in depth by Pitman who has obtained a
number of formulas for distributions related to such Poisson--Kingman
partitions; see Section 8 in \cite{PiPK} or Section 4.5 in~\cite{PiSF}. In particular Pitman has established
sampling formulas in terms of Hermite functions
which provide extensions of the celebrated one due to Ewens \cite{Ewens}.

\subsection{One-type siblings}\label{sec6.2}
We consider now an example related to the fragmentation process at nodes
of the stable tree which has been considered by Miermont~\cite
{Miermont}; see also \cite{AD}. Specifically, we suppose henceforth
that $\beta<2$, and for each fixed value of the parameter $n\in\N$,
all the children of an individual are homebodies
with a probability that decays exponentially in the size of the number
of children, and all children
are migrant otherwise. More precisely, conditionally on $\xi_k=\ell$,
the event $\xi^\mathrm{h}_{n,k}=\ell$ and $\xi^\mathrm{m}_{n,k}=0$ occurs
with probability $\e^{-\ell/n}$, while the event $\xi^\mathrm{h}_{n,k}=0$
and $\xi^\mathrm{m}_{n,k}=\ell$ occurs with probability
$1-\e^{-\ell/n}$.

In this situation, it is easy to check that \eqref{E3} holds with
\[
X^\mathrm{m}_t=\sum_{0<s\leq t} \mathbh{1}_{\{\Delta
X_s>\epsilon
_s\}}
\Delta X_s\quad  \mbox{and}\quad
X^\mathrm{h}_t=X_t-X^\mathrm{m}_t,
\]
where $\Delta X_s$ stands for the size of the jump (if any) of the
$\operatorname{stable}(\beta)$ process~$X$ at time $s$ and $\epsilon_s$ for an
independent standard exponential mark which is attached to each jump of
$X$. Well-known properties of L\'evy processes imply that the
processes~$X^\mathrm{h}$ and $X^\mathrm{m}$ are independent and that
$X^\mathrm{h}$
is an (Esscher) exponential transform of $X$. More precisely, the
bivariate Laplace exponent of $(X^\mathrm{h},X^\mathrm{m})$ is then
given by
\[
\Psi(q,r)=b\bigl((q+1)^{\beta}-1\bigr)+b\bigl(r^{\beta}+1-(r+1)^{\beta}\bigr).
\]

Since the law of $X^\mathrm{h}_t$ is simply an exponential transform of
that of $X_t$,
\[
\P(X^\mathrm{h}_t\in\mathrm{d } x)= t^{-1/\beta}\e^{-tx-bt}\rho
(t^{-1/\beta
}x)\, \mathrm{d } x,
\]
and we deduce from the Ballot theorem that
\[
\Lambda^{(1)}( \mathrm{d } t)= \rho(0) t^{-1-1/\beta} \e
^{-bt}\,\mathrm{d } t.
\]
As the first passage process $T_{\cdot}$ is a subordinator which is
independent of $X^\mathrm{m}$, and $Y=X^\mathrm{m}\circ T$
results from Bochner's subordination, and we get that the L\'evy measure
$\Lambda$ of $(T,Y)$ can be expressed in the form
\[
\Lambda(\mathrm{d}x_1\, \mathrm{d} x_2)=
\rho(0) x_1^{-1-1/\beta} \e^{-bx_1}\,\mathrm{d }x_1\,\P(X^\mathrm
{m}_{x_1}\in
\mathrm{d}x_2).
\]

\subsection{Migration forced by cut-off}\label{sec6.3}
In the preceding two examples, the subordinator $X^\mathrm{m}$
was either deterministic (a pure drift) or independent of the L\'evy
process~$X^\mathrm{h}$.
Our last example shows that more general situations may arise.
Specifically, we consider the parameter $n$
as a threshold and decide that at most $n$ of
the children of each individual are homebodies and the rest are migrants.
In other words, $\xi^\mathrm{h}_{n,k}=\xi_k\wedge n$ and $\xi
^\mathrm{m}_{n,k}=(\xi_k-n)^+$.

Then \eqref{E3} is fulfilled with
\[
X^\mathrm{m}_t=\sum_{0<s\leq t} \mathbh{1}_{\{\Delta X_s>1\}}
(\Delta X_s-1) \quad \mbox{and}\quad
X^\mathrm{h}_t=X_t-X^\mathrm{m}_t.
\]
We stress that the jump times of $X^\mathrm{m}$ are exactly the times when
$X^\mathrm{h}$
has a jump of size 1; in particular the processes $X^\mathrm{h}$ and
$X^\mathrm{m}$ are not independent.

\section*{Acknowledgment} I would like to thank the referee for
having carefully checked the manuscript.

%

%

\printaddresses


\begin{thebibliography}{20}

\bibitem{AD}
%
\begin{barticle}[mr]
\bauthor{\bsnm{Abraham},~\bfnm{Romain}\binits{R.}} \AND
\bauthor{\bsnm{Delmas},~\bfnm{Jean-Fran{\c{c}}ois}\binits{J.-F.}}
(\byear{2008}).
\btitle{Fragmentation associated with {L}\'evy processes using snake}.
\bjournal{Probab. Theory Related Fields}
\bvolume{141}
\bpages{113--154}.
\bid{doi={10.1007/s00440-007-0081-2}, mr={2372967}}
\end{barticle}
%
\endbibitem

\bibitem{ADV}
%
\begin{bmisc}[auto:SpringerTagBib|2009-01-14|16:51:27]
\bauthor{\bsnm{Abraham},~\bfnm{R.}\binits{R.}},
\bauthor{\bsnm{Delmas},~\bfnm{J.~F.}\binits{J.~F.}}
\AND
\bauthor{\bsnm{Voisin},~\bfnm{G.}\binits{G.}}
(\byear{2010}).
\bhowpublished{Pruning a L\'evy continuum random tree. Preprint.
Available at \url{http://hal.archives-ouvertes.fr/hal-00270803\_v2/}}.
\end{bmisc}
%
\endbibitem

\bibitem{LP}
%
\begin{bbook}[mr]
\bauthor{\bsnm{Bertoin},~\bfnm{Jean}\binits{J.}}
(\byear{1996}).
\btitle{L\'evy Processes}.
\bseries{Cambridge Tracts in Mathematics}
\bvolume{121}.
\bpublisher{Cambridge Univ. Press}, \baddress{Cambridge}.
\bid{mr={1406564}}
\end{bbook}
%
\endbibitem

\bibitem{Be1}
%
\begin{barticle}[mr]
\bauthor{\bsnm{Bertoin},~\bfnm{Jean}\binits{J.}}
(\byear{2009}).
\btitle{The structure of the allelic partition of the total population for
{G}alton--{W}atson processes with neutral mutations}.
\bjournal{Ann. Probab.}
\bvolume{37}
\bpages{1502--1523}.
\bid{doi={10.1214/08-AOP441}, mr={2546753}}
\end{barticle}
%
\endbibitem

\bibitem{Be2}
%
\begin{barticle}[mr]
\bauthor{\bsnm{Bertoin},~\bfnm{Jean}\binits{J.}}
(\byear{2010}).
\btitle{A limit theorem for trees of alleles in branching processes
with rare
neutral mutations}.
\bjournal{Stochastic Process. Appl.}
\bvolume{120}
\bpages{678--697}.
\bid{doi={10.1016/j.spa.2010.01.017}, mr={2603059}}
\end{barticle}
%
\endbibitem

\bibitem{Chauvin}
%
\begin{barticle}[mr]
\bauthor{\bsnm{Chauvin},~\bfnm{B.}\binits{B.}}
(\byear{1986}).
\btitle{Sur la propri\'et\'e de branchement}.
\bjournal{Ann. Inst. H. Poincar\'e Probab. Statist.}
\bvolume{22}
\bpages{233--236}.
\bid{mr={850758}}
\end{barticle}
%
\endbibitem

\bibitem{DuLG}
%
\begin{barticle}[mr]
\bauthor{\bsnm{Duquesne},~\bfnm{Thomas}\binits{T.}} \AND
\bauthor{\bsnm{Le~Gall},~\bfnm{Jean-Fran{\c{c}}ois}\binits{J.-F.}}
(\byear{2002}).
\btitle{Random trees, {L}\'evy processes and spatial branching processes}.
\bjournal{Ast\'erisque}
\bvolume{281}
\bpages{1--147}.
\bid{mr={1954248}}
\end{barticle}
%
\endbibitem

\bibitem{Dwass}
%
\begin{barticle}[mr]
\bauthor{\bsnm{Dwass},~\bfnm{Meyer}\binits{M.}}
(\byear{1969}).
\btitle{The total progeny in a branching process and a related random walk.}
\bjournal{J.~Appl. Probab.}
\bvolume{6}
\bpages{682--686}.
\bid{mr={0253433}}
\end{barticle}
%
\endbibitem

\bibitem{Ewens}
%
\begin{barticle}[mr]
\bauthor{\bsnm{Ewens},~\bfnm{W.~J.}\binits{W.~J.}}
(\byear{1972}).
\btitle{The sampling theory of selectively neutral alleles}.
\bjournal{Theoret. Population Biology}
\bvolume{3}
\bpages{87--112}.
\bid{mr={0325177}}
\end{barticle}
%
\endbibitem

\bibitem{GP}
%
\begin{barticle}[mr]
\bauthor{\bsnm{Griffiths},~\bfnm{R.~C.}\binits{R.~C.}} \AND
\bauthor{\bsnm{Pakes},~\bfnm{Anthony~G.}\binits{A.~G.}}
(\byear{1988}).
\btitle{An infinite-alleles version of the simple branching process}.
\bjournal{Adv. in Appl. Probab.}
\bvolume{20}
\bpages{489--524}.
\bid{doi={10.2307/1427033}, mr={955502}}
\end{barticle}
%
\endbibitem

\bibitem{Harris}
%
\begin{bbook}[mr]
\bauthor{\bsnm{Harris},~\bfnm{Theodore~E.}\binits{T.~E.}}
(\byear{1963}).
\btitle{The Theory of Branching Processes}.
\bseries{Die Grundlehren der Mathematischen Wissenschaften, Bd.}
\bvolume{119}.
\bpublisher{Springer}, \baddress{Berlin}.
\bid{mr={0163361}}
\end{bbook}
%
\endbibitem

\bibitem{H}
%
\begin{barticle}[mr]
\bauthor{\bsnm{Hutzenthaler},~\bfnm{Martin}\binits{M.}}
(\byear{2009}).
\btitle{The virgin island model}.
\bjournal{Electron. J. Probab.}
\bvolume{14}
\bpages{1117--1161}.
\bid{mr={2511279}}
\end{barticle}
%
\endbibitem

\bibitem{HW}
%
\begin{barticle}[mr]
\bauthor{\bsnm{Hutzenthaler},~\bfnm{M.}\binits{M.}} \AND
\bauthor{\bsnm{Wakolbinger},~\bfnm{A.}\binits{A.}}
(\byear{2007}).
\btitle{Ergodic behavior of locally regulated branching populations}.
\bjournal{Ann. Appl. Probab.}
\bvolume{17}
\bpages{474--501}.
\bid{doi={10.1214/105051606000000745}, mr={2308333}}
\end{barticle}
%
\endbibitem

\bibitem{Jirina}
%
\begin{barticle}[mr]
\bauthor{\bsnm{Ji{\v{r}}ina},~\bfnm{Miloslav}\binits{M.}}
(\byear{1958}).
\btitle{Stochastic branching processes with continuous state space}.
\bjournal{Czechoslovak Math. J.}
\bvolume{8}
\bpages{292--313}.
\bid{mr={0101554}}
\end{barticle}
%
\endbibitem

\bibitem{Kallenberg}
%
\begin{bbook}[mr]
\bauthor{\bsnm{Kallenberg},~\bfnm{Olav}\binits{O.}}
(\byear{2002}).
\btitle{Foundations of Modern Probability},
\bedition{2nd} ed.
\bpublisher{Springer}, \baddress{New York}.
\bid{mr={1876169}}
\end{bbook}
%
\endbibitem

\bibitem{Miermont}
%
\begin{barticle}[mr]
\bauthor{\bsnm{Miermont},~\bfnm{Gr{\'e}gory}\binits{G.}}
(\byear{2005}).
\btitle{Self-similar fragmentations derived from the stable tree. {II}.
{S}plitting at nodes}.
\bjournal{Probab. Theory Related Fields}
\bvolume{131}
\bpages{341--375}.
\bid{doi={10.1007/s00440-004-0373-8}, mr={2123249}}
\end{barticle}
%
\endbibitem

\bibitem{PiPK}
%
\begin{bincollection}[mr]
\bauthor{\bsnm{Pitman},~\bfnm{Jim}\binits{J.}}
(\byear{2003}).
\btitle{Poisson--{K}ingman partitions}.
In \bbooktitle{Statistics and Science: A {F}estschrift for {T}erry {S}peed}.
\bseries{Institute of Mathematical Statistics Lecture
Notes---Monograph Series}
\bvolume{40}
\bpages{1--34}.
\bpublisher{IMS}, \baddress{Beachwood, OH}.
\bid{doi={10.1214/lnms/1215091133}, mr={2004330}}
\end{bincollection}
%
\endbibitem

\bibitem{PiSF}
%
\begin{bbook}[mr]
\bauthor{\bsnm{Pitman},~\bfnm{J.}\binits{J.}}
(\byear{2006}).
\btitle{Combinatorial Stochastic Processes}.
\bseries{Lecture Notes in Math.}
\bvolume{1875}.
\bpublisher{Springer}, \baddress{Berlin}.
\bid{mr={2245368}}
\end{bbook}
%
\endbibitem

\bibitem{Taib}
%
\begin{bbook}[mr]
\bauthor{\bsnm{Ta{\"{\i}}b},~\bfnm{Ziad}\binits{Z.}}
(\byear{1992}).
\btitle{Branching Processes and Neutral Evolution}.
\bseries{Lecture Notes in Biomathematics}
\bvolume{93}.
\bpublisher{Springer}, \baddress{Berlin}.
\bid{mr={1176317}}
\end{bbook}
%
\endbibitem

\bibitem{Takacs}
%
\begin{bbook}[mr]
\bauthor{\bsnm{Tak{\'a}cs},~\bfnm{Lajos}\binits{L.}}
(\byear{1967}).
\btitle{Combinatorial Methods in the Theory of Stochastic Processes}.
\bpublisher{Wiley}, \baddress{New York}.
\bid{mr={0217858}}
\end{bbook}
%
\endbibitem

\end{thebibliography}
\end{document}